\documentstyle[12pt,amscd,amssymb,righttag]{amsart}
\newcounter{nomer}
\newtheorem{thm}{Theorem}[nomer]
\newtheorem{corol}[thm]{Corollary}

\newtheorem{lemma}[thm]{Lemma}

\theoremstyle{definition}
\newtheorem{defin}[thm]{Definition}

\theoremstyle{remark}
\newtheorem{absatz}[thm]{}

\font\smc=cmcsc10 at 12pt

\def\R {{\Bbb R}}

\def\N{{\Bbb N}}
\def\e{{\epsilon}}
\def\Z {{\Bbb Z}}
\def\s{{\Bbb S}}
\def\Aut{{\mathrm Aut}\,}
\def\Homeo{{\mathrm Homeo}\,}
\def\H{{\mathcal H}}
\def\QED{\nobreak\quad\ifmmode\roman{Q.E.D.}\else{\rm Q.E.D.}\fi}
\def\oskip{\par\vbox to4mm{}\par}
\begin{document}
\noindent {Th\'eorie des groupes/{\it Group theory}}\vskip .5cm
\noindent
\noindent
{\Large\bf
Some extremely amenable groups}
\oskip
Thierry {\smc Giordano} and Vladimir {\smc Pestov}
\oskip

{\it T.G.:} Department of Mathematics and Statistics,
University of Ottawa, Ottawa, Ontario, K1N 6N5, Canada.

E-mail: {\tt giordano@@science.uottawa.ca}

URL: {\tt http://www.science.uottawa.ca/mathstat/proff/giordano.htm}

\oskip

{\it V.P.: }
School of Mathematical and Computing Sciences, Victoria University of
Wellington, P.O. Box 600, Wellington, New Zealand.

E-mail: {\tt vova@@mcs.vuw.ac.nz}

URL:  {\tt http://www.mcs.vuw.ac.nz/$^\sim$vova}\\

\footnotetext{ Nom de la personne qui doit corriger les \'epreuves:
Thierry Giordano \\
Addresse: D\'epartement de math\'ematiques et statistique,
Universit\'e d'Ottawa, \\585 King Edward,
Ottawa, Ontario, K1N 6N5, Canada\\
N$^\circ$ de t\'el\'ephone: +1-613-562-5800 ext. 3514.
N$^\circ$ de t\'el\'ecopieur: +1-613-562-5776.\\
Courier \'electronique: giordano@@science.uottawa.ca}
\smallskip

\oskip

\par
{\small\bf Abstract} {\small\sl --- A topological group $G$ is extremely
amenable
if every continuous action of $G$ on a compact space has a fixed point.
Using the concentration of measure techniques
developed by Gromov and Milman, we prove that
the group of automorphisms of a Lebesgue space with a non-atomic
measure is extremely amenable
with the weak topology but not with the uniform one.
Strengthening a de la Harpe's result, we show that a
von Neumann algebra is approximately finite-dimensional
if and only if its unitary group with the strong
topology is the product of an extremely amenable group with a compact group.}

\oskip
\oskip

\noindent
{\large\bf Quelques groupes extr\^emement moyennables}
\oskip

{\small\bf R\'esum\'e} {\small\it --- Un groupe topologique $G$ est
extr\^emement
moyennable
si toute action continue de $G$ sur un espace compact poss\`ede un point fixe.
En utilisant les techniques de concentration de mesure d\'evelopp\'ees
par Gromov et Milman, nous d\'emontrons que le groupe des automorphismes d'un
espace de Lebesgue avec une mesure diffuse
est extr\^emement moyennable s'il est muni de la topologie
faible, mais ne l'est pas avec la topologie uniforme.
Si $M$ est une alg\`ebre de von Neumann,
nous montrons en utilisant un r\'esultat de P. de la Harpe que $M$ est
approximativement de dimension finie si et seulement si son groupe unitaire
(muni de la topologie forte) est le produit d'un groupe compact et
d'un groupe extr\^emement moyennable.
}

\oskip\medskip


\subsection*{Version fran\c caise abr\'eg\'ee} ---
Un groupe topologique dont toute action continue sur un espace compact
poss\`ede un point fixe est dit {\it extr\^emement moyennable} \cite{Gra}.
Un tel groupe non-trivial
est moyennable, mais n'est jamais localement compact \cite{V}.

Alors que les premiers exemples de tels groupes
ont \'et\'e difficiles \`a
trouver\cite{HC,Ba},
il a \'et\'e r\'ecemment d\'emontr\'e que de nombreux groupes
``de dimensions infinies'' sont extr\^emement moyennables
\cite{GrM,Gl,P1,P5}. Les preuves de ces
r\'esultats utilisent souvent les
techniques de concentration de mesure dans des espaces
de grandes dimensions.

Rappelons qu'un mm{\it -espace} $(X,d,\mu)$ est la
donn\'ee d'un espace m\'etrique
$(X,d)$ muni d'une mesure de probabilit\'e $\mu$.
Une suite g\'en\'eralis\'ee $(X_\alpha,d_\alpha,\mu_\alpha)$
de $mm$-espaces est dite {\it de
L\'evy} si pout tout $\epsilon >0$ et toute famille d'ensembles bor\'eliens
$A_\alpha \subset X_\alpha $ avec $\liminf\mu_\alpha(A_\alpha)>0$, on a
$\lim\mu_\alpha((A_\alpha)_\e)=1$, o\`u $A_\e$ d\'enote
le $\epsilon$-voisinage de $A$.
\\[1.5mm]
{\it Exemples.}
\begin{enumerate}
\item Les groupes de permutations $S_n$ munis de la distance
de Hamming et de la mesure uniforme (Maurey \cite{Ma}). \par
\item Les groupes orthogonaux $SO(n)$ ou les groupes unitaires $SU(n)$
avec la distance induite par la norme de Hilbert-Schmidt et
la mesure de Haar normalis\'ee (Gromov et Milman \cite{GrM}).
\end{enumerate}

Un groupe topologique m\'etrisable est un {\it groupe de L\'evy}
s'il existe
une suite g\'en\'era\-li\-s\'ee $(G_\alpha)_{\alpha \in A}$ de sous-groupes
compacts de $G$ telle que
\begin{enumerate}
\item La suite
$(G, d, \mu_\alpha)$ est de L\'evy, o\`u
$d$ d\'enote une distance
 invariante \`a droite qui induit la topologie de $G$ et o\`u $\mu_\alpha$
 la mesure de Haar normalis\'ee sur $G_\alpha$.
\item Pour tout ensemble fini $g_1, g_2, \cdots, g_N$
de $G$ et tout $\epsilon >0$,
 il existe $\alpha \in A$ tel que, pour tout $\beta \ge \alpha $,
 $d(g_i,G_\beta) < \epsilon$ pour $ 1 \le i \le N.$
\end{enumerate}

G\'en\'eralisant l\'eg\`erement le r\'esultat de \cite{GrM} (cf. \cite{Gl}
et \cite{P5}),
nous avons
\\[1.5mm]
\noindent
{\bf Th\'eor\`eme.}  --- {\it  Un groupe de L\'evy est extr\^emement
moyennable.}
\\[1.5mm]
Soit $(X,\mu)$ un espace bor\'elien standard muni d'une mesure
de probabilit\'e diffuse.
Notons $\Aut^\ast(X,\mu)$ (resp. $\Aut(X,\mu)$)
le groupe des automorphismes
mesurables non-singuliers pr\'eservant la classe de la mesure
$\mu$ (resp. pr\'eservant
$\mu$), muni de la topologie forte.
\\[1.5mm]
\noindent
{\bf Th\'eor\`eme.}  --- {\it  Tant $\Aut^\ast(X,\mu)$ que $\Aut(X,\mu)$
sont des
groupes de L\'evy et donc sont extr\^emement moyennables.}
\\[1.5mm]
\noindent
{\bf Corollaire.}  --- {\it  Si $\mu$ est une mesure $\sigma$-finie,
$\Aut(X,\mu)$ est extr\^emement moyennable.}
\\[1.5mm]
Par contre, $\Aut(X,\mu)$ muni de la topologie uniforme n'est pas
extr\^emement moyennable,
que $\mu$ soit une mesure non-atomique finie ou infinie.

Pour une alg\`ebre de von Neumann $M$, notons $U(M)_s$ son groupe unitaire
muni de la topologie $\sigma (M,M_*)$. Nous avons alors les deux r\'esultats
suivants.
\\[1.5mm]
\noindent
{\bf Th\'eor\`eme.}  --- {\it  Soit $M$ une alg\`ebre de
von Neumann sans partie discr\`ete
finie. Alors $M$ est approximativement de dimension finie si et seulement
si $U(M)_s$
est extr\^emement moyennable.}
\\[1.5mm]
\noindent
{\bf Corollaire.}  --- {\it  Une alg\`ebre de von Neumann est
approximativement de dimension finie si et seulement si son
groupe unitaire $U(M)_s$ est le produit direct
d'un groupe compact et d'un groupe extr\^emement moyennable.}
\\[1.5mm]

Un groupe topologique $G$ est {\it fortement moyennable}
si toute action continue proximale de $G$ sur un espace 
compact poss\`ede un point fixe \cite{Gl1}.
Utilisant un argument semblable 
\`a celui de 
\cite{Pat} (Thm. 2), nous obtenons
\\[1.5mm]
\noindent 
{\bf Corollaire.}  --- {\it  Une C*-alg\`ebre $A$ est nucl\'eaire si 
et seulement si son groupe unitaire $U(A)$, muni de la 
topologie $\sigma (A,A^*)$ (induite) faible,
est fortement moyennable.}

\bigskip

\noindent$\underline{\phantom{xxxxxxxxxxxxxxxxxxxxxxxxxxxxx}}$
\bigskip

\subsection*{\S \arabic{nomer}.  Introduction}
\medskip
Following \cite{Gra} and \cite{Gl}, we call a topological group
whose each continuous action on a compact Hausdorff space has a fixed point,
{\it extremely amenable} or having the {\it fixed point on compacta property.}
Since such a
group satisfies Day's fixed point property, it is amenable.
At the same time, an extremely amenable group is never locally compact, by
this result of Veech \cite{V} (cf. also \cite{Pym}): every locally compact
group $G$ acts freely on a suitable compact space. The first examples were
difficult to find \cite{HC}, \cite{Ba}. However, recently many
`infinite-dimensional' groups were shown to be
extremely amenable: the unitary group $U(\H)_s$ of an infinite-dimensional
Hilbert space with the strong operator topology (Gromov and Milman \cite{GrM}),
the group $L_1(X,U(1))$ of measurable maps from a non-atomic
standard Borel space to the circle group with
the topology of convergence in measure (Glasner \cite{Gl} and
Furstenberg and Weiss, unpublished), the
orientation-preserving
homeomorphism groups $\Homeo_+({\mathbb I})$ and
$\Homeo_+({\mathbb R})$ \cite{P1}, the group of isometries of the Urysohn
metric
space \cite{P5}. Most of the proofs use the technique of concentration of
measure on high-dimensional structures.

If $X$ is a standard Borel space endowed with a measure $\mu$,
we denote by $\Aut(X,\mu)$ (resp. $\Aut^\ast(X,\mu)$)
the group of all measure (resp. measure class) preserving
transformations of $(X,\mu)$ equipped with the weak topology.
In this Note we show that, if $\mu$ is a non-atomic sigma-finite
(resp. finite)
measure, the groups $\Aut(X,\mu)$ and $\Aut^\ast(X,\mu)$
are extremely amenable. With the uniform topology the former
group is shown to be non-amenable.
Strengthening a de la Harpe's result \cite{dlH2},
we show that a
von Neumann algebra $M$ is approximately finite-dimensional (AFD)
if and only if its unitary group $U(M)_s$ with the strong
topology is isomorphic to the product of an extremely amenable group with
a compact one.
As a corollary, we get a new characterization of nuclear $C^\ast$-algebras.

\stepcounter{nomer}
\subsection*{\S \arabic{nomer}. L\'evy families and groups.}
\begin{defin}[\cite{Gr,M2}]
An $mm${\it -space} is a triple $(X,d,\mu)$, where $d$ is a
metric on a set $X$ and $\mu$ is a probability measure on $(X,d)$.
A net $(X_\alpha,d_\alpha,\mu_\alpha)$ of $mm$-spaces
forms a {\it L\'evy family} if,
whenever $A_\alpha\subseteq X_\alpha$ are Borel subsets with
$\liminf\mu_\alpha(A_\alpha)>0$, for each $\e>0$
$\lim\mu_\alpha((A_\alpha)_\e)=1$. (Here $A_\e$ denotes
the $\e$-neighbourhood of a set $A$.)
\end{defin}

\begin{absatz} {\it Examples of L\'evy families.}
1. (Maurey \cite{Ma}.)
The permutation groups $S_n$ of rank $n\geq 1$, equipped
with the uniform measure and the Hamming distance,
\begin{equation}
d_n(\sigma,\tau)=\frac 1n\left\vert\{i\colon \sigma(i)\neq\tau(i).\}
\right\vert
\label{hamming}
\end{equation}
 \par
2. (Gromov and Milman \cite{GrM}.)
The special orthogonal groups $SO(n)$
(or the special unitary groups $SU(n)$), $n\geq 1$, with the
Hilbert-Schmidt operator metric and the normalized Haar measure.
\end{absatz}

The following definition and result slightly generalize those due to
Gromov and Milman \cite{GrM} (cf. \cite{Gl} and \cite{P5}).

\begin{defin}
A metrizable topological group $G$ is a {\it L\'evy group} if there is a net
$(G_\alpha)_{\alpha\in A}$
of compact subgroups of $G$ with the following properties:
\begin{enumerate}
\item The family of $mm$-spaces $(G,d,\mu_\alpha)$ is L\'evy, where
$d$ is any compatible right-invariant metric on $G$ and
$\mu_\alpha$ denotes the normalized Haar measure on $G_\alpha$.
\item For every finite collection $g_1,g_2,\dots,g_N\in G$, $N\in\N$, and
every $\e>0$ there is an $\alpha\in A$ with the property: for all
$\beta\geq\alpha$ and $i=1,2,\dots,N$, $d(g_i,G_\beta)<\e$.
\end{enumerate}
\end{defin}

\begin{thm} Every L\'evy group is extremely amenable.
\qed
\label{gm}
\end{thm}

\stepcounter{nomer}
\subsection*{\S \arabic{nomer}.  Automorphism groups of a Lebesgue space}
Let $X$ be a standard Borel space and let $\mu$ be a non-atomic probability
measure on $X$.
To every non-singular Borel automorphism $T$ of $(X,\mu)$ one associates
a linear isometry $\Phi(T)$ of $L^2(X,\mu)$,
\begin{equation}
\Phi(T)f(x)=\left(\frac{d\mu\circ T^{-1}}{d\mu}(x)\right)^{\frac 12}f(T^{-1}x),
\mbox{ for }x\in X.
\label{quasi}
\end{equation}
The strong operator topology on the unitary group of $L^2(X,\mu)$ induces a
Polish group topology on $\Aut^\ast(X,\mu)$ (and also on $\Aut(X,\mu)$), called
the {\it strong topology}. A finer group topology,
called the {\it uniform topology,}
is determined on both groups by the metric
\[d_{u}(T,S)=\mu\left(\{x\in X\colon T(x)\neq S(x)\}\right).\]

\begin{thm}
\label{main}
The group $\Aut(X,\mu)$ of all measure-preserving
automorphisms of a non-atomic standard Borel
space $X$, equipped with the strong topology, is a L\'evy group.
\end{thm}

\begin{pf} We identify $(X,\mu)$ with $([0,1],\lambda)$, where $\lambda$
is the Lebesgue measure, and $S_{2^n}$
with the subgroup of measure-preserving
automorphisms of $[0,1]$ mapping each dyadic interval of
rank $n$ onto a dyadic interval of rank $n$
via a translation. According to the Weak Approximation Theorem
(cf. e.g. \cite{Hal}, pp. 65--68), the union of the increasing sequence of
such groups is everywhere dense in $\Aut([0,1],\lambda)$.
The restriction of the uniform metric $d_{u}$ to
each $S_{2^n}$ is the Hamming distance (\ref{hamming}).
Now the proof follows from Maurey's result (1.2.(1)) and
Theorem \ref{gm}.
\end{pf}

\begin{corol}
The group $\Aut(X)$ of all measure-preserving
automorphisms of a standard sigma-finite measure
space $(X,\mu)$, equipped with the strong topology,
is extremely amenable.
\end{corol}

\begin{pf}
For $\mu$ finite the result follows from Theorems \ref{main} and \ref{gm},
while for $\mu$ sigma-finite one approximates the group $\Aut(X,\mu)$ with
groups of finite measure-preserving automorphisms and uses a standard
compactness argument.
\end{pf}

Denote by $\Aut(X,\mu)_u$ the group $\Aut(X,\mu)$ equipped with the uniform
topology.

\begin{thm}
The group $\Aut(X,\mu)_u$ is not amenable,
where $\mu$ is either a finite or a sigma-finite non-atomic measure.
\label{notamen}
\end{thm}

\begin{pf}
Let $\mu$ be an invariant probability
measure on $X={\mathrm{SL}}(3,\R)/{\mathrm{SL}}(3,\Z)$.
For a Borel set $A\subseteq X$ of measure $\frac 12$, the function
$h_A=\chi_A-\chi_{A^c}$ belongs to the unit sphere $\s_0$ of
$L^2_0(X,\mu)=\{f\in L^2(X,\mu)\colon \int_X f d\mu=0\}$.
There is a $G$-equivariant positive linear operator of norm $1$
from the space ${\mathrm{UCB}}(\s_0)$ of
uniformly continuous bounded functions on $\s_0$ to the space
of right uniformly continuous bounded functions on
$\Aut(X,\mu)$ (with the left action of the group), under which
a $\varphi\in{\mathrm{UCB}}(\s_0)$ goes to the function
$\Aut(X,\mu)\ni T\mapsto \varphi\left(\Phi(T)h_A \right)$.
(Cf. formula (\ref{quasi}).)
If $\Aut(X,\mu)_u$ were amenable, there would exist an $\Aut(X,\mu)$-invariant
mean on ${\mathrm{UCB}}(\s_0)$, which
would imply (Prop. 3.1 in \cite{P4})
that the representation of ${\mathrm{SL}}(3,\R)$ in $L^2_0(X,\mu)$
is amenable in the sense of Bekka,
and by \cite{B}, Remark 5.10, would have a non-zero invariant
vector, a contradiction.
The infinite case is settled by an argument similar to one in \cite{dlH}.
\end{pf}

Replacing the Weak Approximation Theorem by its generalization based on
a result of Tulcea \cite{IT}, we obtain the following result.

\begin{thm}
\label{main2}
The group $\Aut^\ast(X,\mu)$ of all measure class preserving
automorphisms of a non-atomic standard Borel
space $(X,\mu)$, equipped with the strong topology, is a L\'evy group
and therefore extremely amenable. \qed
\end{thm}


\stepcounter{nomer}
\subsection*{\S \arabic{nomer}. Unitary groups of approximately finite
dimensional von Neumann algebras}

If $M$ is a von Neumann algebra, we denote by $U(M)_s$ its unitary group
endowed with the $\sigma(M,M_\ast)$-topology. If $M$ is acting on a Hilbert
space $\H$, this topology on $U(M)$ coincides with the strong operator
topology.

\begin{thm} Let $M$ be a von Neumann algebra without finite
atomic part. Then $M$ is AFD if and only if its unitary group
$U(M)_s$ is extremely amenable.
\end{thm}

\begin{pf} Sufficiency follows from de la Harpe's
result \cite{dlH2}, in which the separability assumption is not essential
\cite{Haa}. Necessity is a direct consequence of the following
four Lemmas.
\end{pf}

\begin{lemma} The direct product of a family of extremely amenable
groups, equipped with the product topology, is extremely
amenable.
\end{lemma}

Using the structure of type I von Neumann algebras and \cite{P5}, Thm. 2.2,
we have

\begin{lemma} If $M$ is a finite non-atomic type I von Neumann algebra, then
$U(M)_s$ is extremely amenable.
\end{lemma}

\begin{lemma} If $M$ is a finite continuous AFD von Neumann algebra, then
$U(M)_s$ is extremely amenable.
\end{lemma}

\begin{pf} Follows from a generalization of a result by
Glasner--Furstenberg--Weiss (\cite{P5}, Thm. 2.2).
\end{pf}

\begin{lemma} Let $M$ be a properly infinite AFD von Neumann algebra.
Then $U(M)_s$ is extremely amenable.
\end{lemma}

\begin{pf} By Elliott's results \cite{El2}, we can assume that $M$ has a
separable predual and is approximated by an increasing sequence of
finite-dimensional factors. It follows that $U(M)_s$ is a L\'evy group.
\end{pf}

\begin{corol} A von Neumann algebra $M$ is approximately finite dimensional
if and only if its unitary group, $U(M)_s$,
endowed with the strong operator topology, is
the product of a compact group and an extremely amenable group.
\end{corol}

A topological group $G$ is called {\it strongly amenable} if
every continuous proximal action of $G$ on a compact space has a fixed point
\cite{Gl1}. 
By using an argument similar to that in \cite{Pat} (Thm. 2), we obtain
the following.

\begin{corol} A $C^\ast$-algebra $A$ is nuclear if and only if  
its unitary group $U(A)$, equipped with the
topology $\sigma (A,A^\ast)$, is strongly amenable. 
\end{corol} 

%

{\small
\subsection*{Acknowledgements}
The authors gratefully acknowledge
support from the Swiss National Science Foundation and
express their gratitude to Professor Pierre
de la Harpe for his hospitality at
the Universit\'e de Gen\`eve.
T.G. was also partially supported by an NSREC operating grant,
as well as was V.P. during a visit to the University of Ottawa.
The research of V.P. was also supported by the VUW
Research Development Fund and by a Marsden Fund grant
of the Royal Society of New Zealand.

\enddocument
\bye